\input psfig
\def\endofproof{ \hfill $\#$ }
\def\wlog{ without loss of generality }
\centerline{\bf Hubbard Forests}
\bigskip
\bigskip
\centerline{Alfredo Poirier}
\centerline{Mathematics Department}
\centerline{SUNY at StonyBrook}
\centerline{StonyBrook, NY 11794}
\bigskip
\bigskip
{\bf Abstract.} 
{\sl The theory of Hubbard trees (as made precise in [P1]) 
provides an effective classification of non-linear post-critically finite polynomial maps from  ${\bf C}$ to itself. 
This note will extend this classification to the case of maps from a finite union of copies of  ${\bf C}$  to itself. 
Maps which are post-critically finite and nowhere linear will be characterized by
 a ``forest'',  
which is made up out of one tree in each copy of ${\bf C}$.}  
\vskip 1in
\centerline{\bf Introduction.}
\bigskip
Let ${\bf C}^\cup$ be the disjoint union of finitely many copies of the complex numbers ${\bf C}$. 
We will consider proper holomorphic maps $f: {\bf C}^\cup \to {\bf C}^\cup$. 
In other words the map $f$ restricted to a copy ${\bf C}_u$ of the complex numbers, 
is a polynomial taking values in the (perhaps different) copy ${\bf C}_{u'}$. 
In order to avoid trivial behavior we impose the condition that the degree of each of these polynomials is at least $2$ (compare also remark 0.6). 
For such functions the concepts of critical and post-critical points are well defined, 
as well as the notions of filled Julia set $K(f)$, Julia set $J(f)$ and Fatou components (compare [M]). 
For example the {\it filled Julia set $K(f)$} is by definition the set of points 
which remain bounded under iteration, 
and the {\it Julia set} is then the boundary of the filled Julia set $K(P)$. 
The following are elementary results in this direction. 
In both cases the proof is based on some elementary topology and in the fact 
that all such sets are dynamically defined in terms of {\it the basin of attraction of `$\infty$'}. 
\bigskip
{\bf 0.1. Lemma.}
{\it The Julia set and the filled Julia set are completely invariant under the map $f$.} 
\endofproof
\bigskip
{\bf 0.2. Lemma.}
{\it The map $f$ carries each Fatou component onto a Fatou component as a branched covering map.}
\endofproof
\bigskip
We are concerned in this note only with the postcritically finite case. 
An effective classification of this kind of map can be established in analogy with that given for postcritically finite polynomials by Hubbard Trees 
as follows (compare [DH], [P1], [P2] and [P3]). 
Let $\Omega(f)$ be the set of critical points of $f$. 
In this case the {\it orbit} ${\cal O}(\Omega(f))$ of the critical set is by definition finite. 
For every copy ${\bf C}_u$ of the complex plane included in ${\bf C}^\cup$, 
denote by {\bf Post}$(u)$ the portion of this orbit ${\cal O}(\Omega(f))$ that belongs to ${\bf C}_u$. 
We can define the {\it regulated tree} $T_u$ as the regulated convex closure of the set {\bf Post}$(u)$ in $K(f)\cap {\bf C}_u$ (compare [DH] or [P1], see also Appendix B). 
Thus we have a natural generalization of Hubbard Trees. 
(Compare Theorem A.) 
It follows easily from lemma s 0.1 and 0.2 that $f(T_u) \subset T_{u'}$ where 
${\bf C}_{u'}$ is the range of the polynomial map $f|_{{\bf C}_u}$. 
We are also interested in the abstract analogue of the above construction. 
\bigskip
By a {\it graph} $G$ is meant a finite simplicial complex of dimension $\le 1$. 
(In particular,  $G$  is compact, the set of vertices  $V$  is a finite subset
  of  $G$ , and  $G-V$  is a finite union of open 1-cells called {\it edges}.)
\medskip
{\bf 0.3. Definition.} 
An {\it abstract Hubbard Forest} will mean a graph $H$, 
together with a function $\phi$ from the set of vertices $V \subset H$ to itself, 
a {\it degree function} $d(v)>0$ from $V$ to the positive integers, 
and an {\it angle function}  $\angle(e,e') \in {\bf R/Z}$  between pairs of edges incident 
at a common vertex, 
satisfying the following six conditions. 
(These are all immediate 
generalizations of the corresponding conditions when  $H$  is connected.)
\bigskip\noindent
(1) Each connected component of $H$ is simply connected (i.e, a point or a tree).
\smallskip\noindent
(2) Each connected component contains at least one `critical vertex', 
that is a vertex with $d(v)>1$.
\smallskip\noindent
(3) Two adjacent vertices map to distinct vertices belonging to a
common connected component.
\smallskip\noindent
(4) {\it Angle Condition:} 
The angle $\angle(e,e')$ should be skew-symmetric, with 
$\angle(e,e')=0$ if and only if $e=e'$, and with 
$\angle_v(e,e'')=\angle_v(e,e')+\angle_v(e',e'')$ 
whenever three edges are incident at a vertex $v$. 
As usual we extend $\phi$ to a map from $H$ to itself whose restriction to each edge is a homeomorphism. 
We require also this angle to be compatible with the `dynamics' $\phi$ and the degree 
$d$ in the following way, 
$$
\angle_{\phi(v)}(\phi(e),\phi(e'))=d(v)\angle_v(e,e') \hbox{\it \hskip 0.5in (mod 1).}
$$
(In this case $\phi(e),\phi(e')$ are incident at a vertex $\phi(v)$ where the angle is measured.) 
\smallskip\noindent 
(5) {\it Normalization at periodic Julia vertices:} 
We say that a vertex $v \in V$ is a {\it Julia vertex} if  
there are no periodic critical vertices in its forward orbit. 
We require then that 
for any periodic Julia vertex $v$ where $m$ edges meet, 
the angle between two edges incident at $v$ 
should be a non trivial multiple of $1/m$. 
(A vertex which is not Julia will be called a {\it Fatou vertex}.)
\smallskip\noindent
(6) {\it Expanding Condition:} 
We define the {\it distance $d_{\bf H}(v,v')$} between vertices 
as the number of edges between $v$ and $v'$ or as $\infty$ if $v$ and $v'$ belong to different components of $H$. 
We require then that 
for any pair of different periodic Julia vertices $v$ and $v'$ 
there is a $k \ge 0$ such that 
$d_{\bf H}(\phi^{\circ k}(v),\phi^{\circ k}(v'))>1$. 
\bigskip
By an {\it isomorphism} between abstract Hubbard Forests is meant a homeomorphism from one graph to the other which preserves all of this structure. 
(In particular, it must carry any vertex of degree  $d$ to a vertex of degree  $d$, compare $\S$2.2.)
\bigskip
\eject
{\bf  Theorem A.}  
{\it Let  $f$  be a nowhere linear postcritically finite proper  holomorphic map from a finite union  ${\bf C}^\cup$ of copies of  ${\bf C}$  to itself. 
Let  $M$  be a finite forward invariant set containing all critical points.
Then the regulated hull  H  of  M  can be canonically given the structure
  of an abstract Hubbard forest. 
The set  $V$  of vertices of this forest will consist of  M  together with finitely many points where three or  more edges come together.} 
\endofproof
\bigskip
\bigskip
{\bf Theorem B.} 
{\it Conversely, every abstract Hubbard forest is isomorphic to one which is constructed in this way, 
for some map $f : {\bf C}^\cup \to {\bf C}^\cup$  and some finite subset of  $K(f)$. 
Here ${\bf C}^\cup$  consists of one copy of  ${\bf C}$  for each component of  H . 
Furthermore, this map  f  is unique up to component-wise affine conjugation.} 
\bigskip
Theorem A follows from our previous discussion. 
For the missing details the reader should take a look at Section 3 and Remark B.7. 
Theorem B is a reformulation of Theorem 2.4. 
Clearly the key step is to prove that every cycle of Trees can be realized; 
and then proceed backwards. 
\bigskip
\bigskip
In section 1 we will introduce some notation. 
We will also prove a result about realization of coverings, 
which will allow us to inductively proceed backwards. 
This last result is also used in sections 2 and 4. 
In section 2 we prove (assuming the results from section 5) 
the {\it Realization Theorem for Abstract Hubbard Forests}. 
In section 3 we recall some standard results about Hubbard Trees; proofs may be found either in [P1] or [P3]. 
In section 4 we establish auxiliary results to be used in section 5. 
In section 5 we prove that a cyclic Hubbard Forest can be realized. 
For this we need the following result due to Peter Jones. 
\bigskip
{\bf 0.4. Theorem ([J, Corollary 2].} 
{\it Let $P$ be a Postcritically Finite Polynomial with Julia set $J(P)$. 
Suppose $\phi:{\bf C} \to {\bf C}$ is a homeomorphism holomorphic off $J(P)$.\break
Then $\phi$ is an affine map.}
\bigskip 
For the convenience of the reader we include two appendices. 
In Appendix A we follow Milnor for the basic notions of mapping schemata for the study of polynomial maps as described above (compare [M, $\S$2]). 
In Appendix B we follow Douady and Hubbard for the study of regulated sets in compact, connected and locally connected 
subsets of the complex plane (compare [DH, Chapter 2]). 
\bigskip
{\bf 0.5 Sketch of the Proof.} 
The idea of the proof of Theorem B is as follows. 
We assume that the given Hubbard Forest is cyclic. 
In this case it is easy to construct `a first return map' from a `canonical extension' of a connected component $H_u$ of $H$ to itself 
(compare $\S$4.1 and $\S$4.4). 
It follows from our construction that this is an expanding Hubbard Tree and therefore there is a unique polynomial $P_u$ associated to this component. 
Now, if this construction is made simultaneously at all connected components, 
it follows that 
the associated polynomials $P_{u_1},\dots P_{u_r}$ all have the same degree (compare $\S$4.1). 
This set of polynomials will represent the `first return map' of the postcritically finite map of\break 
Theorem B.  
Therefore to complete the proof we need a `factorization lemma'. 
This is done with the help of Theorem 0.4 using the filled Julia set $K(P_u)$ of $P_u$.  
Suppose the connected component $H_{u_1}$ of $H$ maps to the connected component $H_{u_2}$. 
Following the `patterns' of $K(P_{u_1})$ and $K(P_{u_2})$ 
(which `realize' $H_{u_1}$ and $H_{u_2}$) 
we construct a continuous branched covering map from the copy $C_{u_1}$ of the complex plane 
(on which the polynomial $P_{u_1}$ is defined) 
to the copy $C_{u_2}$, 
which is analytic off $J(P_{u_1})$ and maps Fatou components onto Fatou components as a branched covering map  
(compare Proposition 5.1). 
We use Theorem 0.4 to show that this map is a polynomial of the appropriate degree 
and the result follows (compare $\S$4.6-10). 
\bigskip
{\bf 0.6 Remark.} 
It should be noted that condition (2) above can be relaxed. 
In fact, we need only assume that every cycle of trees contains a critical point. 
The proofs that follow, apply without any important modification. 
\bigskip
\bigskip
{\bf Acknowledgement.} 
We will like to thank John Milnor for helpful discussions. 
Also, we want to thank the Geometry Center, University of Minnesota
and Universidad Cat\'olica del Per\'u 
for their material support. 
\bigskip
\bigskip
\centerline{\bf 1. Abstract Coverings.}
\bigskip
\bigskip 
{\bf 1.1. Angled Trees.} 
By an {\it (angled) tree H=(T,V,\angle)} will be meant a finite connected acyclic $m$-dimensional simplicial complex where $m=0$ or $m=1$, 
together with a function 
$\ell,\ell' \mapsto \angle(\ell,\ell')=\angle_v(\ell,\ell') \in {\bf Q/Z}$ 
which assigns a rational modulo 1 to each pair of edges $\ell,\ell'$ which meet at a common vertex $v$. 
This angle $\angle(\ell,\ell')$ should be skew-symmetric, with $\angle(\ell,\ell')=0$ if and only if $\ell=\ell'$, 
and with $\angle_v(\ell,\ell'')=\angle_v(\ell,\ell')+\angle_v(\ell',\ell'')$ whenever three edges are incident at a vertex $v$. 
Here $V$ represents the set of vertices and $T$ the underlying topological tree 
(or a point if $m=0$). 
In what follows the case $m=0$ is easy and special, 
and will be excluded from most of the discussion. 
\medskip
Such an angle function determines a preferred isotopy class of embeddings of $T$ into ${\bf C}$, 
in which the induced cyclic order of edges at a vertex $v$ is preserved. 
More precisely, the set $E_v$ of edges incident at a vertex $v \in V$ has a natural cyclic order determined by the angle function. 
{\bf Definition.} 
{\it An admissible embedding of $T$} 
is an embedding at which the cyclic order (taken counterclockwise) of the edges incident at a common vertex induced by the embedding 
is the same as that induced by the angle function. 
\medskip
Once an admissible embedding is given we can walk around the tree in the counterclockwise sense. 
Of course the order in which the vertices will be found is independent of the admissible embedding. 
Formally, we have the following. 
{\bf Definition.} 
If $\ell,\ell' \in E_v$ are consecutive in the cyclic order of $E_v$, 
we say that $(v,\ell,\ell')$ is a {\it pseudoaccess to $v$}. 
(Note that $\ell=\ell'$ if there is only one incident edge.) 
Take a pseudoaccess $(v,\ell,\ell')$ to $v$, 
and let the end points of the edge $\ell'$ be $v,v' \in V$. 
At $E_{v'}$ let $\ell''$ be the successor of $\ell'$ in the cyclic order. 
We say that the pseudoaccess $(v',\ell',\ell'')$ is the {\it successor} of $(v,\ell,\ell')$. 
\bigskip
{\bf 1.2 Lemma.}
{\it Let $(T,V,\angle)$ be an angled tree. 
The successor function in the set of pseudoaccesses to the vertices in $V$ is a complete cyclic order. 
Furthermore if $|V|$ denotes the number of vertices in $V$, the number of pseudoaccesses is $2(|V|-1)$. }
\medskip
{\bf Proof.}
A trivial induction on the cardinality of $V$. 
\endofproof
\bigskip
\bigskip
{\bf 1.3. Definition.} 
Let $H_i=(T_i,V_i,\angle_i)$ $(i=1,2)$ be angled trees. 
By an\break
{\it abstract $H_1$-covering of $H_2$} will be meant a function 
$f:V_1 \to V_2$, together with 
a local degree 
$\delta_f:V_1 \to {\bf Z}$ which assigns a positive integer $\delta_f(v)\ge 1$ to every vertex $v \in V_1$. 
These functions are subject to the following restrictions. 
\smallskip
Whenever $v$ and $v'$ are the endpoints of a common edge $\ell$ in $H_1$, 
we require that\break
$f(v) \ne f(v')$. 
Next, extend $f$ to a map $f:H_1 \to H_2$ which carries each edge homeomorphically onto the shortest path joining the images of its endpoints. 
We require then that 
$\angle_{f(v)}(f(\ell),f(\ell'))=\delta_f(v)\angle_v(\ell,\ell')$. 
(In this case $f(\ell)$ and $f(\ell')$ are incident at the vertex $f(v)$ where the angle is measured.) 
\smallskip
By definition a vertex {\it $v \in V_1$ is critical} if $\delta(v)>1$ and 
{\it non-critical} otherwise. 
The number $deg(f)=1+\sum_{v \in V_1}(\delta(v)-1)$ is the {\it degree of the covering $f$}. 
\medskip
A particular case is when $deg(f)=1$. 
In this case we say that $H_2$ is an {\it extension of $H_1$}. 
Here we will think of $H_1$ as a subset of $H_2$. 
(Compare Lemma 1.5 where this abuse of notation is used). 
\bigskip
{\bf 1.4. Homogeneous Coverings.} 
Unfortunately, coverings do not necessarily map consecutive edges in $E_v$ to consecutive edges in $E_{f(v)}$.
Therefore we can not expect in general to have an induced map between pseudoaccesses. 
Furthermore, even if consecutive edges map to consecutive edges 
the induced map between pseudoaccesses is not necessarily onto. 
{\bf Definition.} 
We say that an abstract $H_1$-covering of $H_2$ is {\it homogeneous} 
if consecutive edges in $E_v$ map to consecutive edges in $E_{f(v)}$ and 
the induced map between pseudoaccesses is a $deg(f)$ to 1 order preserving covering. 
An alternative definition (compare Corollary 1.7) is that every vertex in $V_2$ has $n=deg(f)$ inverses, 
counting multiplicity by $\delta_f$. 
In general it is easier to work with homogeneous coverings. 
This will represent no loss of generality as the following lemma shows. 
\bigskip
{\bf 1.5. Lemma.} 
{\it Let $(f,\delta_f)$ be an abstract 
$H_1$-covering of $H_2$ of degree $n>0$. 
Then there exists an extension angled tree $H_1'$ of $H_1$ and 
$(f',\delta_{f'})$ a homogeneous  
abstract $H'_1$-covering of $H_2$ also of degree n, 
such that 
$f(v)=f'(v)$ and $\delta_f(v)=\delta_{f'}(v)$ for every $v \in V_1$. 
Furthermore such extension is unique.
} 
\medskip
{\bf Proof.} 
Without loss of generality after subdividing $H_1$ we can assume that if $v,v'$ are the end points of an edge in $H_1$ then $f(v),f(v')$ are also the end points of an edge in $H_2$. 
We will construct in a canonical way an extension of $H_1$ and an associated abstract covering. 
The induced map will define an order preserving covering between pseudoaccesses. 
Topological reasons will allow us to deduce that in fact, the degree should be $n$. 
\smallskip
We start by constructing an extension of $H_1$. 
For this we take any pseudoaccess $(v,\ell,\ell')$ in $H_1$. 
Then by hypothesis $f(\ell) \in E_{f(v)}$. 
Let ${\cal E}$ be the successor of $f(\ell)$ in $E_{f(v)}$ and 
$\theta=\angle_{f(v)}(f(\ell),{\cal E})$ be the angle between these two edges. 
Then by definition of abstract covering $\angle_{v}(\ell,\ell')\ge \theta/\delta(v)$. 
If equality holds, we make no modification and proceed to the successor of 
$(v,\ell,\ell')$ in the order of pseudoaccesses. 
If equality does not hold we take the connected component of $T_2-\{f(v)\}$ which contains ${\cal E}$ 
and `glue' it at $v$ forming an angle of $\theta/\delta(v)$ with $\ell$. 
With the obvious definitions this is clearly an abstract covering of $H_2$. 
Note that no new edges or vertices would be added at this new `branch' at $v$ if we perform a similar construction. 
After walking around the tree following that order of the pseudoaccesses of $H_1$ we finish with an extension $H'_1$ of $H_1$ and of the covering for which the 
induced map between pseudoaccesses is an order preserving $m$ to $1$ covering. 
We need only to prove that $n=m$. 
\smallskip
Now let $|V_1|,|V_2|$ be the number of vertices in $V'_1$ and $V_2$ respectively. 
Then Lemma 1.2 tell us that $2(|V_1|-1)=2m(|V_2|-1)$. 
Also note that a pseudoaccesses at $v$ of angle $\theta_v$ is mapped to a pseudoaccess at $f(v)$ of angle $\delta(v)\theta_v$. 
Thus adding this relation for all accesses in the tree $H'_1$ we get 
$|V_1|+n-1=m|V_2|$. 
Now this relation and the former give $n=m$ and the result follows. 
\endofproof
\bigskip
{\bf 1.6 Corollary.} 
{\it Let $(f,\delta_f)$ be an abstract $H_1$-covering of $H_2$ of degree $n>0$. 
Then every point $v \in V_2$ has at most $n$ inverses, counting multiplicity by $\delta_f$.} 
\medskip
{\bf Proof.} 
If not, there is no way to construct a homogeneous $H_1'$-covering of $H_2$. 
\endofproof
\bigskip
{\bf 1.7 Corollary.} 
{\it Let $(f,\delta_f)$ be an abstract $H_1$-covering of $H_2$ of degree $n>0$. 
Then $(f,\delta_f)$ is homogeneous if and only if 
for every $v \in V_2$ 
$$
n=\sum_{v' \in V_1:f(v')=v} \delta_f(v'). 
$$
In other words every vertex in $V_2$ has as many inverses (counting `multiplicity') as the degree of $f$.} 
\medskip
{\bf Proof.} 
No new vertices can be added in the construction of an extension. 
\endofproof
\bigskip
{\bf 1.8 Proposition.} 
{\it Let $(f,\delta_f)$ be a homogeneous abstract $H_1$-covering of $H_2$ of degree $n>0$. 
For any admissible embedding of $T_2$ in the complex plane, 
there is an admissible embedding of $T_1$ 
and a polynomial $P$ of degree n which realizes $(f,\delta_f)$. 
Furthermore, 
such polynomial is unique up to precomposition with an affine map.}
\medskip
{\bf Proof.} 
Without loss of generality we assume that the edges of the embedded tree $T_2$ are analytical arcs. 
Choose an end of $T_2$ 
(i.e, a vertex at which a unique edge terminates). 
We join this vertex to $\infty$ by an analytic arc without crossing the tree, 
and cut open this configuration. 
Then glue together $n$ copies of this surface following the pattern of $T_1$. 
Clearly with the induced complex structure, this is a complex plane. 
Furthermore, the induced identification map is a degree $n$ proper holomorphic map and therefore a polynomial $P$. 
By construction the local degree of $P$ agrees with $\delta_f$ at embedded vertices. 
As we only used necessary conditions for this construction, 
we have that after normalization $P$ is unique. 
\endofproof
\bigskip
We will also need the following technical lemma about composition of abstract\break
coverings. 
\medskip
{\bf 1.9 Lemma.} 
{\it Let $H_1,H_2,H_3$ be abstract angled trees. 
Suppose $(f_1,\delta_1)$ is a\break
homogeneous abstract $H_1$-covering of $H_2$, and 
$(f_2,\delta_2)$ an abstract $H_2$-covering of $H_3$.\break
Define $g=f_2 \circ f_1:V_1 \to V_3$, and 
$\delta(v)=\delta_1(v)\delta_2(f_1(v))$ for $v \in V_1$. 
Then $(g,\delta)$ is an abstract $H_1$-covering of $H_3$ 
of degree $deg(f_1)deg(f_2)$.} 
\medskip
We will denote such abstract coverings as $(f_2,\delta_2)*(f_1,\delta_1)$. 
Note that $*$ is associative whenever this makes sense. 
\medskip
{\bf Proof.} 
Clearly $(g,\delta)$ is a $H_1$-abstract covering of $H_3$. 
We must verify that\break
$deg(g)=deg(f_1)deg(f_2)$. 
For this we recall first that by Corollary 1.7, 
for every $v' \in V_2$ we have 
$deg_1(f)=\sum_{v \in V_1:f_1(v)=v'} \delta_1(v)$. Then 
$$
\eqalign{deg(g)-1 &= \sum_{v \in V_1}(\delta(v)-1) \cr
&= \sum_{v \in V_1}(\delta_2(f_1(v))\delta_1(v)-1) \cr
&= \sum_{v \in V_1}(\delta_2(f_1(v))-1)\delta_1(v)+ \sum_{v \in V_1}(\delta_1(v)-1) \cr
&= \sum_{v' \in V_2} \sum_{\{v \in V_1:f_1(v)=v'\}}(\delta_2(v')-1)\delta_1(v)+ deg(f_1)-1 \cr
&= \sum_{v' \in V_2} (\delta_2(v')-1) \sum_{\{v \in V_1:f_1(v)=v'\}}\delta_1(v)+ deg(f_1)-1 \cr
&= \sum_{v' \in V_2} (\delta_2(v')-1)deg(f_1) + deg(f_1)-1 \cr
&= (deg(f_2)-1)deg(f_1) + deg(f_1)-1 \cr\cr
&= deg(f_2)deg(f_1)-1.}
$$
\endofproof 
\bigskip
\bigskip
\centerline{\bf 2. Hubbard Forests.} 
\bigskip
\bigskip
{\bf 2.1. Definition.} 
By an {\it abstract Hubbard Forest ${\bf H^*}$ 
with {\it underlying mapping schema} $S({\bf H^*})=(|S|,F,w)$}
will be meant 
a collection $H_u$ ($u \in |S|$) of angled trees and 
a collection $(f_u,\delta_u)$ of abstract $H_u$-coverings of $H_{F(u)}$ of degree $d(u)=w(u)+1$; 
satisfying the following (usual) conditions. 
\smallskip
Let $V=\bigcup_{u \in |S|} V_u$. 
We have a well defined function $f:V \to V$ between vertices 
as defined by $f(v)=f_u(v)$ whenever $v \in V_u$. 
We say that a vertex $v \in V$ is a {\it Julia vertex} if  
there are no periodic critical vertices in its forward orbit. 
We require then that 
for any periodic Julia vertex $v$ where $m$ edges meet, 
the angle between two edges incident at $v$ 
should be a non trivial multiple of $1/m$. 
(A vertex which is not Julia will be called a {\it Fatou vertex}.)
\smallskip
Let $v,v'$ be vertices in the same angled tree $H_u$, 
we define the {\it distance $d_{\bf H^*}(v,v')$} between these vertices 
as the number of edges between $v$ and $v'$. 
We require then that 
for any pair of different periodic Julia vertices $v$ and $v'$ 
belonging to the same tree $H_u$ there is a $k \ge 0$ such that 
$d_{\bf H^*}(f^{\circ k}(v),f^{\circ k}(v'))>1$. 
\bigskip
{\bf 2.2 Definition: Equivalent Hubbard Forests.}
Let ${\bf H^*,{H^*}'}$ be two abstract Hubbard Forests. 
We say that they are {\it equivalent} (in symbols ${\bf H^*} \cong {\bf {H^*}'}$) if 
the following (obvious) conditions are satisfied: 

(1) There exists a bijection 
$\phi:|S| \to |S'|$ which establishes an equivalence between 
the underlying mapping schemata $S({\bf H^*})=(|S|,F,w)$ and 
$S'({\bf {H^*}'})=(|S'|,F',w')$ 
(compare Appendix A). 

(2) For each $u \in |S|$ there is a homeomorphism $\phi_u:T_u \to T'_{\phi(u)}$ 
which preserves 

  (i) vertices, i.e, $\phi_u(V_u) = V'_{\phi(u)}$, 

 (ii) dynamics, i.e, $f_{\phi(u)}(\phi_u(v))=\phi_{F(u)}(f_u(v))$, 

(iii) degree, i.e, $\delta_u(v)=\delta'_{\phi(u)}(\phi_u(v))$ for all $v \in V_u$, and 

(iv) angles, i.e, $\angle_v(\ell,\ell')=\angle_{\phi_u(v)}(\phi_u(\ell),\phi_u(\ell'))$ for all edges $\ell,\ell' \in E_v$ (Here $v \in V_u$). 
\bigskip
{\bf 2.3. The space ${\cal P}^S$ of polynomial maps.} 
Let $S=(|S|,F,w)$ be a reduced mapping schema 
(see Appendix A for definitions). 
Form the disjoint union $|S|\times {\bf C}$ of $n$ copies of the complex numbers ${\bf C}$, 
where $n$ is the number of vertices of $S$. 
In other words, 
replace each vertex of $S$ by a copy of ${\bf C}$. 
Let ${\cal P}^S$ be the space consisting of all maps $f$ from $|S|\times {\bf C}$ to itself 
such that the restriction of $f$ to each component $u \times {\bf C}$ is a monic centered polynomial of degree $d(u)=w(u)+1$, taking values in $F(u) \times {\bf C}$. 
Clearly the spaces ${\cal P}^S$ and ${\cal P}^{S'}$ are isomorphic if and only if 
the defining mapping schemata are equivalent (in the sense described in Appendix A). 
\medskip
{\bf Remark.}
Theorem A says that whenever $f \in {\cal P}^S$ is postcritically finite, 
then any invariant set $M$ containing the critical set $\Omega(f)$ naturally defines an abstract Hubbard Forest ${\bf H^*}_{f,M}$. 
\bigskip
{\bf 2.4. Theorem.} 
{\it An abstract Hubbard Forest ${\bf H^*}$ 
(with underlying mapping schema $S({\bf H^*})=(|S|,F,w))$, 
can be realized by a Postcritically Finite $f \in {\cal P}^{S'}$ 
and an invariant set $M$ containing the critical set $\Omega(f)$
(i.e, ${\bf H^*} \cong {\bf H^*}_{f,M}$) 
if and only if the reduced mapping schemata $S({\bf H^*})$ and ${S'}$ are equivalent. 
Furthermore, once a bijection $\phi:|S| \to |S'|$ which realizes the equivalence has been chosen 
then $f$ is unique modulo affine change of coordinates in each copy ${\bf C}_{u'}$ 
(here ${u'} \in |S'|$).}
\bigskip
As was pointed out before, 
the main technical problem is the realization of cycles of Trees. 
For this we need the following Theorem whose proof we delay until section 5. 
\bigskip
\eject
{\bf 2.5. Theorem.} 
{\it A cyclic abstract Hubbard Forest ${\bf H^*}$ 
(with underlying mapping schema $S({\bf H^*})=(|S|,F,w))$, 
can be realized by a Postcritically Finite $f \in {\cal P}^{S'}$ 
and an invariant set $M$ containing the critical set $\Omega(f)$
(i.e, ${\bf H^*} \cong {\bf H^*}_{f,M}$) 
if and only if the reduced mapping schemata $S({\bf H^*})$ and ${S'}$ are equivalent. 
Furthermore, once a bijection $\phi:|S| \to |S'|$ which realizes the equivalence has been chosen 
then $f$ is unique modulo affine change of coordinates in each copy ${\bf C}_{u'}$ 
(here ${u'} \in |S'|$).}
\bigskip
{\bf Proof of Theorem 2.4.} 
Clearly we may realize each component of $S({\bf H^*})$\break 
independently. 
For this we realize first the cyclic parts of $S({\bf H^*})$ using Theorem 2.5. 
Next, we inductively proceed backwards using Lemma 1.5 and Proposition 1.8. 
The result\break 
follows. 
\endofproof
\bigskip
\bigskip
\centerline{\bf 3. Some results on Hubbard Trees.} 
\bigskip
\bigskip
{\bf 3.1.} 
Whenever the underlying mapping schema $S$ consists of a single vertex, 
the definition of abstract Hubbard Forest 
coincides with that of abstract Hubbard Tree. 
In this case we know that Hubbard Trees can be realized. 
In what follows we let ${\bf H}=(T,V,\angle)$ be the Hubbard Tree with dynamics $f$ and local degree function $\delta$. 
\medskip
{\bf Theorem.} (See [P1].)  
{\it Let ${\bf H}$ be an abstract Hubbard tree. 
Then there is a postcritically finite polynomial P and an invariant set $M \supset \Omega(P)$ such that ${\bf H}_{P,M} \cong {\bf H}$. 
Furthermore, P is unique up to affine conjugation.} 
\endofproof
\bigskip
This abstract Hubbard tree also gives information about external rays as the following Theorem essentially due to Douady and Hubbard shows 
(see [DH, Chap VII] or [P1]).
\bigskip
{\bf 3.2. Theorem.}
{\it The number of rays which land at a periodic Julia vertex $v$ is equal to 
the number of incident edges of the tree $T$ at $v$, 
and in fact, 
there is exactly one ray landing between each pair of consecutive edges. 
Furthermore, the ray which lands at $v$ between $\ell$ and $\ell'$ maps to the ray 
which lands at $f(v)$ between $f(\ell)$ and $f(\ell')$.}
\endofproof
\bigskip
{\bf Remark.} After these Theorems there is no reason to distinguish between the abstract Hubbard Tree and the unique polynomial which realizes it. 
\medskip
{\bf Definition.} 
A point $p \in J(P)$ is {\it terminal} if there is only one external ray landing at $p$. 
Otherwise $p$ is an {\it incidence point}. 
For incidence points we distinguish between {\it branching} (if there are more than two rays landing at $p$) 
and {\it non branching} (exactly two rays landing at $p$). 
For a postcritically finite polynomial $P$, 
every branching point must be periodic or preperiodic. 
Also every periodic branching point is present as a vertex in any tree ${\bf H}_{P,M}$. 
\bigskip
{\bf 3.3. Proposition.} 
{\it Let P be a Postcritically Finite Polynomial and $z \in J(P)$ a\break 
branching point. 
Then $z$ is preperiodic (or periodic).}
\endofproof
\bigskip
\eject 
{\bf 3.4. Proposition.} 
{\it Let P be a Postcritically Finite Polynomial and $z \in J(P)$ a\break 
periodic incidence point. 
Let $M$ be an invariant set containing the critical points of P, 
then $z \in T_{P,M}$. 
Furthermore, the number of components of $T_{P,M}-\{z\}$ is independent of M and equals the number of components of $J(P)-\{z\}$.}
\endofproof
\bigskip
\bigskip
\centerline{\bf 4. Auxiliary Results.} 
\bigskip
\bigskip
We will establish in this section some auxiliary notation and results. 
We assume that ${\bf H^*}$ is a cyclic Hubbard Forest; i.e, 
its underlying schema $S({\bf H^*})$ is cyclic. 
We write this cycle as $u_0 \mapsto u_1 \dots \mapsto u_r=u_0$. 
We omit the trivial case in which each $V_u$ consists of a single vertex.
\bigskip
{\bf 4.1. A canonical collection} $P_u$ ($u \in |S|$) {\bf of Polynomials.} 
By definition we have a collection of abstract coverings $(f_i,\delta_i)$ each of degree $d(u_i)=w(u_i)+1$. 
Denote by $m=d(u_0)\dots d(u_{r-1})$ the product of these degrees. 
This number $m$ is {\it the degree of ${\bf H^*}$}. 
Inductively, using Lemma 1.5 we may assume without loss of generality for 
$i=1, \dots r-1$ that 
$(f_i,\delta_i)$ is a homogeneous $H_i$-covering of $H_{i+1}$. 
Furthermore, 
again by Lemma 1.5 there is an extension $H'_0$ of $H_0$ such that 
$(f_0,\delta_0)$ is a homogeneous $H'_0$-covering of $H_{1}$. 
Furthermore $H'_0$ with that dynamics and degree induced from 
$(f_{r-1},\delta_{r-1})*\dots *(f_0,\delta_0)$ is clearly an expanding Hubbard Tree. 
By Lemma 1.8, the degree of this abstract Hubbard Tree is $m$.  
Therefore to every $u \in |S|$ 
we can canonically associate a postcritically finite polynomial $P_u$ of degree $m$ 
from the copy ${\bf C}_{u}$ of the complex numbers to itself
(unique up to affine conjugation). 
First we use this fact in order to prove some extension properties for ${\bf H}^*$ 
and afterwards we will proceed to normalize these polynomials. 
\bigskip
{\bf 4.2. Proposition.} 
{\it Let $v \in J(P_u)$ be a periodic point for $P_u:{\bf C}_u \to {\bf C}_u$. 
Then, there is an extension of ${\bf H^*}$ 
on which an orbit corresponding to $v$ is present.}
\medskip
In other words, we may assume that there are vertices $v \in V_u$, $f(v) \in V_{F(u)}$ and so on, 
with the required property. 
\bigskip
{\bf Proof.} 
Take $v \in {\bf C}_0$ periodic under $P_0$. 
We assume that $v$ does not belong to the orbit of the critical set $\Omega(P_0)$ 
(otherwise there is nothing to prove). 
By definition, the angled tree $H'_0$ 
(compare $\S$4.1) 
can be realized as the 
regulated convex closure of a collection of vertices in $K(P_0)$ 
(compare Appendix B). 
We extend this regulated set so that it also includes the forward orbit of $v$ under $P_0$. 
Let ${\hat H}_0$ be this angled tree. 
Note that the dynamics and local degree of $P_0$ defines an abstract $H'_0$ covering of ${\hat H}_0$ of degree $m$. 
Note that there are two ways to extend the angled tree $H'_0$ into an homogeneous covering of ${\hat H}_0$ of degree $m$. 
The first is given by the definition of the Hubbard Tree corresponding to $P_0$. 
Clearly this extension must contain ${\hat H}_0$. 
On the other hand by Lemma 1.5 we can produce extensions 
${\hat H}_i$ of $H_i$ such that 
$(f_i,\delta_i)$ is a homogeneous ${\hat H}_i$-covering of ${\hat H}_{i+1}$ for 
$i=1,\dots,r-1$, and an extension 
$H''_0$ of $H'_0$ so that 
$(f_0,\delta_0)$ is a homogeneous ${H}''_0$-covering of ${\hat H}_{1}$. 
Furthermore ${H}''_0$ with that dynamics and degree induced from 
$(f_{r-1},\delta_{r-1})*\dots *(f_0,\delta_0)$ is an homogeneous abstract covering of ${\hat H}_0$ of degree $m$ and 
therefore must coincide with the previous one as guaranteed by the uniqueness part of Lemma 1.5. 
In other words, one of the $m$ preimages of $P_0(v)$ under 
$(f_{r-1},\delta_{r-1})*\dots *(f_0,\delta_0)$ is canonically identified with $v$. 
The result then follows easily. 
That the resulting collection of abstract coverings is in fact a Hubbard Forest is an easy consequence of the angle condition and Proposition 3.4. 
\endofproof
\bigskip
{\bf 4.3 Corollary.} 
{\it Let $v$ be a preperiodic point for $P_u$. 
Then, there is an extension of ${\bf H^*}$ 
(still a Hubbard Forest based on the same schema $S$) 
on which an orbit corresponding to $v$ is present.}
\endofproof
\bigskip
{\bf 4.4. Normalization of} $\{P_u\}_{u \in |S|}$. 
Now we normalize every $P_u$ (as defined in $\S$4.1) so that they become monic and centered polynomials of degree $m$, as follows. 
By definition every underlying topological tree $T_u$ and set of vertices $V_u$ can be realized within the filled Julia set $K(P_u)$. 
Let $p$ be the landing point of the zero ray $R_0$ corresponding to $P_u$. 
Because of Proposition 4.2, 
we may assume that $p \in V_u$. 
Now we make use of Theorem 3.2. 
Corresponding to $H_0$ this ray $R_0$ lands between two edges $\ell,\ell'$ incident at $p$. 
Using the angle condition, 
we have that between the edges $f_0(\ell)$ and $f_0(\ell')$ incident at $f_0(p)$ 
a fixed ray ${\cal R}$ for the polynomial $P_{F(u)}$ must land. 
After affine change of coordinates in ${\bf C}_{F(u)}$ (if necessary) 
we take this fixed ray ${\cal R}$ to be the zero ray for $P_{F(u)}$. 
Finally we proceed forward in order to associate a `compatible' zero ray to every polynomial $P_{u_i}$. 
Thus in order to prove the existence part of Theorem 2.5 
we only need to construct a 
collection of polynomials $f_0,f_1,\dots,f_r=f_0$ each of degree $d(u_i)$ such that 
$$
P_i=f_{i-1} \circ \dots \circ f_1 \circ f_0 \circ f_{r-1} \dots \circ f_i. 
$$
This will be done in Section 5. 
In the meanwhile, with the preceding normalization we have the following result. 
\medskip
{\bf 4.5 Proposition.} 
{\it Suppose that between the edges $\ell,\ell'$ incident at $v \in V_i$ 
the rational ray $R_\theta$ (for $P_i$) lands. 
Then between the edges $f_i(\ell),f_i(\ell')$ incident at $f_i(v) \in V_{i+1}$ 
the rational ray $R_{d(u_i)\theta}$ (for $P_{i+1}$) lands.}
\medskip
{\bf Proof.} 
The proof is analogous to that given in [P1] for Hubbard Trees and therefore we will only sketch it. 
(All coverings are assume here to be homogeneous.) 
Let $v \in V_u$ be a Julia vertex. 
By definition an {\it access} at $v$ is a pseudoaccess $(v,\ell,\ell')$ 
for which a unique rational ray for $P_{u}$ lands. 
In other words the is a unique rational ray landing between the edges $\ell$ and $\ell'$ incident at $v$. 
Note that the angle condition and Theorem 3.2 imply that the induced map from $(f_{u_i},\delta_{u_i})$ 
between pseudoaccess, would map accesses at $v$ to accesses at $f_{u_i}(v)$. 
Now, given an angled tree $H_i$, there is a unique way 
(once the access corresponding to the argument $0$ is chosen)
to assign an argument to every access in such way that the induced map from 
$(f_{i-1},\delta_{i-1})*\dots *(f_0,\delta_0)*\dots*(f_{i},\delta_{i})$
becomes multiplication by $m$. 
As the assignment $d(u_i)\theta$ is compatible with the dynamics of accesses in $H_{i+1}$, 
the result follows. 
\endofproof
\bigskip
\eject
{\bf 4.6.} 
We take now any given covering $(f_u,\delta_u)$ in the definition of cyclic Hubbard Forest. 
There is no loss of generality to assume that this abstract covering is homogeneous. 
Up to this point we have realized the angled trees $H_u$ and $H_{F(u)}$ as regulated convex sets 
in the filled Julia sets $K(P_u)$ and $K(P_{F(u)})$ of $P_u$ and $P_{F(u)}$ respectively. 
In other words the sets of vertices $V_u,V_{F(u)}$ can be thought as collections of points contained within 
$K(P_u)$ and $K(P_{F(u)})$ respectively. 
Similarly, the underlying topological trees $T_u,T_{F(u)}$ become now the trees
generated by these sets of vertices in $K(P_u) \subset {\bf C}_u$ and $K(P_{F(u)}) \subset {\bf C}_{F(u)}$ respectively 
(compare Remark B.7). 
\smallskip
On the other hand, when this realization of $H_{F(u)}$ has been fixed, 
Proposition 1.8 shows how to realize the abstract covering $(f_u,\delta_u)$ 
as a polynomial ${\cal F}:{\bf C}_u \to {\bf C}_{F(u)}$ of degree $d(u)$ in an essentially unique way. 
Therefore there is second way to realize the angled tree $H_u$ in ${\bf C}_u$. 
In other words, 
there is a tree $T \subset {\bf C}_{u}$, 
a set of vertices $V \subset T$ and 
an essentially unique polynomial map ${\cal F}:{\bf C}_{u} \to {\bf C}_{F(u)}$ of degree $d(u)$
which realizes $(f_u,\delta_u)$ (compare Proposition 1.8). 
The main goal is to prove that these two realizations of $H_u$ within ${\bf C}_u$ can be chosen to coincide 
(compare Proposition 5.1). 
For this, first 
we have to establish some notation and results that will be useful when comparing $T$ with $T_u$. 
\smallskip
For simplicity we set $K({\cal F})={\cal F}^{-1}K(P_{F(u)})$ and 
$J({\cal F})=\partial K({\cal F})={\cal F}^{-1}J(P_{F(u)})$. 
Note that ${\cal F}$ restricted to any connected component of ${\bf C}-J({\cal F})$ is a branched covering map onto its image 
(which is necessarily a Fatou component for $P_{F(u)}$). 
Thus each connected component of ${\bf C}-J({\cal F})$ shall be called a 
{\it `Fatou component'} of ${\cal F}$.  
It follows that angles at Fatou points in $V \subset T$ are actually the angles between radial segments in the uniformizing coordinates of components. 
In particular, the polynomial ${\cal F}$ restricted to ${\bf C}-K({\cal F})$ is by definition a $d(u)$ unbranched covering of
${\bf C}-K(P_{F(u)})$ 
(or equivalently the branching takes place at the removable singularity at $\infty$). 
Thus `external rays' for $K({\cal F})$ are just pullbacks of external rays for $K(P_{F(u)})$. 
We should label these rays in such a way that they convey the same information 
as the external rays for the polynomial $P_u$. 
In essence this means that we must carefully chose the assignment of the ray with argument $0$. 
For this we must use the polynomial $P_u$ as reference: 
Let $p$ be the point were the zero ray for $P_u$ lands. 
Assume without loss of generality (compare Proposition 4.2) that $p \in V_u$. 
Note that by definition there is a point $\alpha(p) \in V \subset T$ which canonically corresponds to $p$. 
Clearly $\alpha(p) \in J({\cal F})$. 
As $p$ is a periodic Julia vertex for $P_u$, 
the corresponding $\alpha(p)$ is not a critical point of ${\cal F}$. 
Thus, among the rays landing at $\alpha(p)$ only one maps under ${\cal F}$ to the zero external ray of  $K(P_{F(u)})$. 
By definition this is the `external zero ray' for $K({\cal F})$. 
As a matter of notation, 
the ray with argument $\theta$ for $K({\cal F})$ will be denoted by $\alpha(R_\theta)$. 
\medskip
Of course, the concepts of branching and incidence points in $K({\cal F})$ 
are well defined because this last set is also locally connected. 
For future reference we state the following result about `external rays' for $K({\cal F})$. 
The proofs are easy consequences of the fact that ${\cal F}$ is a branched covering map. 
\bigskip
{\bf 4.7 Lemma.} 
{\it Let $q' \in J({\cal F})$; 
then the number of rays landing at $q'$ and the number of connected components of $J({\cal F})-\{q'\}$ coincide. 
Furthermore, if $q'$ is a branching point of $K({\cal F})$; 
then only rays with rational argument land at $q'$.} 
\endofproof
\bigskip
\eject
{\bf 4.8.} 
Now let $W=\{q \in K(P_u):\hbox{\it $q$ is preperiodic for $P_u$}\}$. 
Analogously define the set 
$\alpha(W)=\{\alpha(q) \in K({\cal F}):\hbox{\it ${\cal F}(\alpha(q))$ is preperiodic for $P_{F(u)}$}\}$.
There is a natural one to one correspondence between points in $W$ and $\alpha(W)$ defined as follows. 
Take $q \in W$, which because of Corollary 4.3 we assume in $V_u$. 
By definition there is a unique point $\alpha(q)\in \alpha(W)$ which corresponds to $q$. 
Conversely if $\alpha(q) \in \alpha(W)$, 
then we may assume that ${\cal F}(\alpha(q)) \in V_{F(u)}$ and the result follows in an analogous way. 
Furthermore, note that this correspondence associates Fatou (respectively Julia) points to Fatou (respectively Julia) points. 
\smallskip 
Using the `Fatou points' in $\alpha(W)$ as the centers of `Fatou components' of $K({\cal F})$, 
we can define regulated paths in $K({\cal F})$ and therefore regulated sets 
(compare Appendix B). 
It follows easily that a regulated set in $K({\cal F})$ maps to a regulated set in $K(P_{F(u)})$. 
As a consequence of this, 
we have that regulated paths $[q_1,q_2]$ joining points $q_1,q_2 \in W$ 
are in one to one correspondence with 
regulated paths $[\alpha(q_1),\alpha(q_2)]$ joining points $\alpha(q_1),\alpha(q_2) \in \alpha(W)$. 
In other words the vertices in $[q_1,q_2]\cap W$ are in precise 
(order preserving) correspondence with 
the vertices in $[\alpha(q_1),\alpha(q_2)] \cap \alpha(W)$ 
(and viceversa). 
\bigskip
{\bf 4.9 Lemma.}
{\it The ray $R_\theta$ lands at $v \in W \cap J(P_u)$ 
if and only if 
the ray $\alpha(R_\theta)$ lands at the vertex $\alpha(v) \in \alpha(W) \cap J({\cal F})$.}
\medskip
{\bf Proof.} 
This is an easy consequence of the definition of homogeneous coverings, the fact that 
${\cal F}$ is a branched covering and Lemma 4.5. 
\endofproof
\bigskip
{\bf 4.10 Lemma.}
{\it The rays $R_\theta$ and $R_{\theta'}$ land at the same point if and only if 
the rays $\alpha(R_\theta)$ and $\alpha(R_{\theta'})$ land at the same point.}
\medskip
{\bf Proof.} 
The case when either $\theta$ or $\theta'$ is rational follows from Lemma 4.9 and the fact that ${\cal F}$ is a covering map. 
So let $\theta \ne \theta'$ be irrational arguments and suppose that the rays $R_\theta$ and $R_{\theta'}$ both land at $p \in J(P_u)$. 
Let $\alpha(q)$ be the landing point of $\alpha(R_\theta)$ and 
$\alpha(q')$ be that of $\alpha(R_\theta)$. 
We will assume $\alpha(q) \ne \alpha(q')$ and derive a contradiction. 
\smallskip
Let $[\alpha(q),\alpha(q')]$ be the regulated path joining $\alpha(q)$ and $\alpha(q')$ within $K({\cal F})$. 
First we show that $[\alpha(q),\alpha(q')]$ does not contain any center $\alpha(v)$ of a `Fatou Component' $\Delta_{\alpha(v)}$ in $K({\cal F})$. 
To prove this, we assume again otherwise and derive a contradiction. 
Note that by definition, $\overline{\Delta_{\alpha(v)}}-[\alpha(q),\alpha(q')]$ consists of exactly two disjoint components. 
In each of these two components 
there are points $\alpha(z_1),\alpha(z_2)$ at which rays 
$\alpha(R_{\gamma_1}),\alpha(R_{\gamma_2})$ with rational argument land. 
The corresponding rays $R_{\gamma_1},R_{\gamma_2}$ 
then land at different points $z_1 \ne z_2$ in the boundary of the corresponding 
Fatou component with center $v$ in $K(P_u)$. 
Thus we have that $\theta,\theta'$ belong to different connected components of 
${\bf R/Z}-\{\gamma_1,\gamma_2\}$. 
But of course this is impossible. 
\smallskip
Thus we are left with the case where $[\alpha(q),\alpha(q')] \subset J({\cal F})$. 
We take any point\break 
$\alpha(q'') \in [\alpha(q),\alpha(q')]$ different from $\alpha(q)$ and $\alpha(q')$. 
In any neighborhood $U$ of $\alpha(q'') \in J({\cal F})$ 
there is a point $\alpha(z)$ at which a ray $\alpha(R_\beta)$ with rational argument lands 
(this follows from local connectivity). 
If $U$ was chosen small enough, 
this implies (again by local connectivity) that either 
this landing point $\alpha(z)$ belongs to $[\alpha(q),\alpha(q')]$, or 
there is a branching point (again $\alpha(z)$) in $[\alpha(q),\alpha(q')]$. 
In either case the rays $\alpha(R_{\theta_1}),\dots,\alpha(R_{\theta_r})$
landing at $\alpha(z) \in [\alpha(q),\alpha(q')]$ have rational argument by Lemma 4.7.
As $\alpha(q),\alpha(q')$ belong to different components of $[\alpha(q),\alpha(q')]-\{\alpha(z)\}$, 
it follows by Lemma 4.7 that $r>1$ and that  
$\theta,\theta'$ belong to different connected components of ${\bf R/Z}-\{\theta_1,\dots,\theta_r\}$. 
Thus, as the rational rays $R_{\theta_1},\dots,R_{\theta_r}$ land at the same point $z \in J(P_u)$ by Lemma 4.9, 
it follows that $R_{\theta}$ and $R_{\theta'}$ land at different points. 
But this last statement is a contradiction. 
The converse is proven in an analogous way. 
\endofproof
\bigskip
\bigskip
\centerline{\bf 5. Realizing Cyclic Hubbard Forests.} 
\bigskip
\bigskip
In this section we use the notation and results established above. 
The key technical result is the following Proposition. 
\bigskip
{\bf 5.1 Proposition.} 
{\it The abstract covering $(f_u,\delta_u)$ can be realized as a polynomial $f$ of degree $d(u)=deg(f_u)$ 
in the following sense, 
\smallskip
i) $f(v)=f_u(v)$ for all $v \in V_u$; 
\smallskip
ii) the local degree of $f$ at $v$ is $\delta_u(v)$; 
\smallskip
iii) $f(T_u) \subset T_{F(u)}$.
\smallskip
Furthermore such $f$ is unique.}
\medskip
{\bf Proof.} 
We assume without loss of generality 
that $(f_u,\delta_u)$ is a homogeneous covering. 
We proceed as follows. 
When we fix the embedding $T_{F(u)}$, 
by Proposition 1.8 there is a tree $T \subset {\bf C}_{u}$, 
a set of vertices $V \subset T$ and 
an (essentially unique) polynomial map ${\cal F}:{\bf C}_{u} \to {\bf C}_{F(u)}$ of degree $d(u)$
which realizes $(f_u,\delta_u)$ (in the sense described above). 
We need to prove that we can assign coordinates to ${\bf C}_{u}$ in such way that $T_u=T$ and $V_u=V$. 
For this we will construct a homeomorphism $\phi$ of ${\bf C}_{u}$ which maps $T_u$ to $T$, $V_u$ to $V$ and 
is holomorphic off $J(P_u)$. 
The result then follows from Theorem 0.4. 
We will construct $\phi$ in two steps. 
First in the closure of the basin of attraction of infinity $\overline{A_u(\infty)}$, 
and then in the filled Julia set $K(P_u)$. 
\smallskip
First take $\phi$ as the analytic homeomorphism between ${\bf \hat C}-K(P_u)$ and
${\bf \hat C}-K({\cal F})$ (with $\phi(\infty)=\infty$), 
which makes the rays with argument $0$ correspond.  
Lemma 4.10 shows that $\phi$ extends as a homeomorphism between 
$\overline{{\bf \hat C}-K(P_u)}$ and
$\overline{{\bf \hat C}-K({\cal F})}$. 
Furthermore, Lemma 4.9 shows that ${\cal F}(\phi(v))=f_u(v)$.  
\medskip
Now let $v_1,\dots,v_r \in V_u$ be the Fatou points in $V_u$. 
We enumerate the centers of Fatou components of $P_u$, 
subject to the condition that 
for $i > r$, there is a $j<i$ such that $P_u(v_i)=v_j$. 
We will inductively define $\phi$ holomorphically in the Fatou component $U_i$ with center $v_i$, 
in such way that it extends continuously to $\partial U_i$ in a compatible way. 
For this (using Corollary 4.3), 
we assume that $v_i \in V_u$. 
Thus, there are corresponding edges $\ell$ and $\alpha(\ell)$ of $T_u$ and $T$ respectively 
which are radial segments in the uniformizing coordinates with centers $v_i$ and $\alpha(v_i)$. 
Making them correspond in uniformizing coordinates, 
we have an analytic homeomorphism $\phi$ between these Fatou components, 
which extends to the boundary. 
Because ${\cal F}$ is a covering of this set $\phi(U_i)$ onto its image in ${\bf C}_{F(u)}$, 
the angle condition shows 
that $\phi$ is independent of this edge $\ell$. 
We still need to verify that for boundary points, 
this definition of $\phi$ coincides with that given before. 
Note that it is enough to verify this for preperiodic points in $\partial U_i$. 
For this we take any preperiodic point $p \in  \partial U_i$ which we assume with loss of generality to be in $V_u$. 
But then by construction, 
$\phi(p)=\alpha(p)$ in either definition.  
Thus, this extension is compatible. 
\smallskip
After defining $\phi$ in all Fatou components, 
$\phi$ is a homeomorphism of ${\bf C}_u$ which is holomorphic off $J(P_u)$ and therefore affine by Theorem 0.4. 
This completes the proof of Proposition 5.1. 
\endofproof

\bigskip{\bf 5.2 Proof of Theorem 2.4.} 
Proposition 5.1 shows that there is a collection\break
$f_0,f_1,\dots,f_r=f_0$ of polynomials each of degree $d(u_i)$ such that 
$$
P_i=f_{i-1} \circ \dots \circ f_1 \circ f_0 \circ f_{r-1} \circ \dots \circ f_i. 
$$
\noindent
According to Milnor (compare [M, Lemma 2.7]), 
they can normalized in such way that $(f_0,f_1,\dots,f_{r-1}) \in {\cal P}^S$; 
furthermore it is also shown there that this normalization is unique modulo 
the required properties. 
\endofproof
\bigskip
\bigskip
\centerline{\bf Appendix A: Mapping Schemata.} 
\bigskip
\bigskip
{\bf A.1 Definition.} 
By a {\it mapping schema} $S=(|S|,F,w)$ we mean: 

(1) a finite set $|S|$ of points, together with 

(2) a function $F$ from $|S|$ to itself, and also 

(3) a ``weight function" $w$ which assigns an integer $w(v) \ge 0$ called the {\it critical weight} 
to each $v \in |S|$. 

Equivalently, 
such a mapping schema can be represented by a finite graph with one vertex for each $v \in |S|$, 
and with exactly one directed edge $e_v$ leading out from each vertex 
$v$ to a vertex $F(v)$. 
By definition the degree associated with the edge $e_v$ 
(or with the vertex $v$) is the integer $d(v)=w(v)+1 \ge 1$. 
The {\it weight $w(S)$ of a schema} $S$ is by definition the number 
$w(S)=\sum_{v \in |S|} w(v)$. 
\medskip
Such a mapping schema is {\it reduced} if every vertex is critical. 
Suppose that we start with a mapping schema $S$ which satisfies the following very mild condition: 
{\it 
Every cycle in S contains at least one critical vertex.} 
Then there is an 
{\it associated reduced mapping schema} $\bar S$ 
which is obtained from $S$ simply by discarding all vertices of weight zero and shrinking every edge of degree one to a point. 
Note that $S$ and $\bar S$ have the same total weight. 
\bigskip
Any hyperbolic polynomial map $f$ from ${\bf C}$ to ${\bf C}$, or from a finite union of copies of ${\bf C}$ to itself of degree two or more in each copy, 
gives rise to an associated 
{\it full mapping schema} $S=(|S|,F,d)$ by the following construction. 
Let $W(f)$ be the union of the basins of attraction for all attraction periodic orbits of f in ${\bf C}$. 
(Equivalently, 
since $f$ is hyperbolic, 
$W(f)$ is the interior of the filled Julia set $K(P)$.) 
Let $W^{PC}$ be the union of those components of $W(f)$ whose intersection with the orbit 
${\cal O}(\Omega(f))$ of the critical set $\Omega(f)$ is non empty. 
\smallskip
{\bf Definition: The full mapping schema of $f$.} 
This schema $S(f)$ has one vertex $v$ corresponding to every component $W_v \subset W^{PC}$. 
The {\it weight} $w(v)$ is defined as the number of critical points in $W_v$ (counting multiplicity). 
Every vertex $v$ is joined to a vertex $F(v)$ by an edge $e_v$ of degree $w(v)+1$, 
where $F(v)$ is the vertex associated with the component 
$f(W_v)-W_{F(v)} \subset W^{PC}$. 
\bigskip
{\bf A.2 Definition: Equivalent schemata.} 
We say that two mapping schemata\break 
$S=(|S|,F,w)$ and $S'=(|S'|,F',w')$ are 
{\it equivalent} 
if  there is a bijection $\phi:|S| \to |S'|$ between the respective vertices which preserves weight and 
conjugates the mappings $F$ and $F'$. 
\bigskip
{\bf A.3 Remark and Definition.} 
The set $Aut(S)$ of all automorphisms of $S$ which establish an equivalence of $S$ to itself, 
carries a natural group structure. 
For more details about this group $Aut(S)$ we refer to [M]. 
\bigskip
\bigskip
\centerline{\bf Appendix B: Regulated Sets.} 
\bigskip
\bigskip
{\bf B.1} 
In this appendix we consider only a special kind of compact subsets of the complex plane. 
Let $K \subset {\bf C}$ be a compact set with the following properties. 
\smallskip
a) $K$ is connected and locally connected. 

b) Its complement ${\bf C}-K$ consists of a unique component (which of course, is unbounded). 
\medskip
For example if $P:{\bf C} \to {\bf C}$ is a postcritically finite polynomial, it is known that its filled Julia set $K(P)$ satisfies these properties. 
\medskip
Now, the interior $int(K)$ of $K$ has at most a countable number of components. 
We take a point $z_i$ in each of these components $U_i$, 
which we call the {\it center of $U_i$}. 
As $K$ is compact, it follows that every $U_i$ is conformally equivalent to the unit disk $D$. 
Now, the uniformizing coordinate (or Riemann Map) $\phi_i:U_i \to D$
can be chosen in such way that $\phi_i(z_i)=0$. 
Furthermore we have that $\phi_i$ is unique up to a rotation in $D$. 
Therefore we can define {\it internal rays} in $U_i$ as `pull backs' of radial segments in $D$. 
It follows that given two points in the closure $\overline{U_i}$ of $U_i$, 
they can be joined in a unique way by a Jordan arc consisting of (at most two) segments of internal rays. 
We call such arcs (following Douady and Hubbard) {\it regulated}. 
\medskip
As $K$ is locally connected in a compact metric space, 
then $K$ is also arcwise connected. 
This means that given two points $z_1,z_2 \in K$ there is a continuous injective map\break 
$\gamma:I=[0,1]\mapsto K$ such that $\gamma(0)=z_1$ and $\gamma(1)=z_2$. 
Such arcs (really their images) can be chosen in a unique way so that the intersection with the closure of a component $U_i$ is regulated
(see [DH, chapter 2]). 
We still call such arcs {\it regulated}, and denote them by $[z_1,z_2]_K$. 
The following immediate properties hold for regulated arcs (compare also [DH, chapter 2]).
\bigskip
{\bf B.2 Lemma.}
{\it Let $\gamma_1,\gamma_2$ be regulated arcs, then $\gamma_1 \cap \gamma_2$ is regulated.}
\endofproof
\bigskip
{\bf B.3 Lemma.}
{\it Every subarc of a regulated arc is regulated.}
\endofproof
\bigskip
{\bf B.4 Lemma.}
{\it Let $z_1,z_2,z_3 \in K$, then there exists $p \in K$ such that 
\break 
$[z_1,z_2]_K \cap [z_2,z_3]_K=[p,z_2]_K$. 
In particular if $[z_1,z_2]_K \cap [z_2,z_3]_K= \{z_2\}$, 
the set 
\break 
$[z_1,z_2]_K \cup [z_2,z_3]_K$ is a regulated arc.}
\endofproof
\bigskip
{\bf B.5 Regulated Sets.} 
We say that a subset $X \subset K$ is {\it regulated connected} 
if for every $z_1,z_2 \in X$ we have $[z_1,z_2]_K \subset X$. 
We define the {\it regulated hull $[X]_K$ of $X \subset K(P)$} as the minimal closed regulated connected subset of $K$ containing $X$.
\bigskip
{\bf B.6 Proposition.}
{\it If $z_1,...,z_n$ are points in K, 
the regulated hull $[z_1,...,z_n]_K$ of $\{z_1,...,z_n\}$ is a finite topological tree.}
\endofproof
\medskip 
The proof this last Proposition is an easy induction argument using condition b). 
\bigskip
{\bf B.7. Remark: Hubbard Trees (following Douady and Hubbard).} 
(Compare [DH] and [P1].) 
Let $P:{\bf C} \to {\bf C}$ be a postcritically finite polynomial. 
In case the filled Julia set $K(P)$ has non empty interior, each of these components is a component of the Fatou set. 
Furthermore, by definition each bounded Fatou component belongs to $K(P)$. 
Each of these Fatou components has a unique distinguish point which eventually maps to a periodic critical point. 
We will pick always this point $z(U)$ as the center of $U$. 
Hubbard Trees were defined by Douady and Hubbard as the regulated convex closure of the orbit ${\cal O}(\Omega(P))$ of the critical set $\Omega(P)$ of $P$. 
The angle function at Fatou points is defined using the uniformizing coordinate. 
At a Julia point $v$ it is defined as a multiple of $1/m$ where $m$ is the number of rays landing at $v$.  
\bigskip
\bigskip
\centerline{\bf References.} 
\bigskip
\bigskip
[DH] A.Douady and J.Hubbard, \'Etude dynamique des Polyn\^omes Complexes, part I; 
Publ Math. Orsay 1984-1985. 
\smallskip
[J] P. Jones, On Removable Sets for Sobolev Spaces in the Plane. 
Preprint $\#$1991/22 IMS, SUNY@StonyBrook.
\smallskip
[M] J.Milnor, Hyperbolic Components in Spaces of Polynomial Maps; Preprint\break 
$\#$1992/3 IMS, SUNY@StonyBrook. 
\smallskip
[P1] A.Poirier, On Postcritically Finite Polynomials; Thesis, SUNY@StonyBrook, 1993 (to appear). 
\smallskip
[P2] A.Poirier, On the Realization of Fixed Point Portraits; Preprint $\#$1991/20 IMS, SUNY@StonyBrook. 
\smallskip
[P3] A.Poirier, Hubbard Trees; {\it ...in preparation}. 
\bigskip
\bigskip
Typeset in $\TeX$ (August 13, 1992). 
\end

\centerline{\bf 6. Tuning.} 
\bigskip
\bigskip
An important application of the theory of Hubbard Forest is given by the construction below. 
This construction is a step towards the generalization of the 
Douady-Hubbard theory  of `modulation' or `tuning' for  quadratic Polynomials. 
\bigskip
{\bf 6.1} 
Let $P$ be a postcritically finite Polynomial. 
Let $|S|=\{v_1,\dots,v_m\}$ be an invariant set of Fatou points. 
We impose the technical assumption that no critical point in 
$\Omega(P)-|S|$ eventually maps into $|S|$. 
With that induced local degree and dynamics from $P$, 
this set $|S|$ can be given the structure of a mapping schema $S$. 
By hypothesis, this mapping schema 
has the property that every cycle contains a critical point. 
Therefore, there is an associated space ${\cal P}^S$. 
\medskip
On the other hand, 
to each Fatou component $U_i$ (for which the `marked point' $v_i$ is the center), 
we can associate ``zero internal rays" ${\cal R}(i)$ 
with the property  that they map to each other following that induced order from $S$. 
Of course, in general there are many ways for the election of this `boundary marking' 
(compare [M]). 
\medskip
Note that every $f \in {\cal P}^S$ 
is by definition a collection of monic centered polynomials $\{f_i\}$ of degree $d(v_i)$ 
from the copy $\{v_i\}\times {\bf C}$ of the complex numbers to the copy 
$\{P(v_i)\}\times {\bf C}$. 
Therefore if $f$ is assumed also to be in the connectedness locus of ${\cal P}^S$, 
there is a canonically well defined `zero external ray' ${\cal R}'_0(i)$ for 
$K_i(f)=K({\cal F}) \cap (\{v_i\}\times {\bf C})$. 
\medskip
{\bf 6.2 Formal Tuning.} 
In what follows, we fix a boundary marking for this postcritically finite Polynomial $P$. 
Also let $f$ be in the connectedness locus of ${\cal P}^S$. 
Suppose further that the Julia set $J(f)$ is locally connected. 
We replace the Fatou component $U_i$ to which $v_i$ belong, 
by a copy of the dynamical plane $\{v_i\} \times {\bf C}$, 
so that internal angles in $U_i$ correspond to external angles of $K_i(f)$. 
Next we shrink these gluing rays to a point. 
When this is done for all $v_i \in |S|$, 
we inductively proceed backwards until all Fatou components that eventually map 
to some $U_i$ are covered, and repeat the same construction. 
By Moore's Theorem the resulting topological space is homeomorphic to ${\bf R}^2$. 
An we have an induced dynamics in this set. 
This is by definition the {\it formal tuning of $P$ by $f$}. 
 The question is of course if we can give a complex structure to this homeomorphic copy of ${\bf R}^2$, 
so that the induced dynamics is a Polynomial. 
\smallskip
Two questions immediately arise. 
What if $J(f)$ is not locally connected? 
How can we distinguish between two answers that are merely topologically conjugate? 
\medskip
{\bf 6.3 Hyperbolic Tuning.} 
Assume now that $P$ is postcritically finite and hyperbolic; 
i.e, every critical point of $P$ eventually maps to a critical cycle. 
Let $H=H(P)$ be the Hubbard Tree of $P$. 
Let $|S|$ be the set of all Fatou points in $H$. 
Then for every $f \in {\cal P}^S$ postcritically finite, 
we will construct a new polynomial $P \bot \hbox{ $f$}$ which would be called 
{\it $P$ tuned by $f$ (respect to the fixed boundary marking)}. 
This is done as follows. 
\smallskip
First we need to extend the Hubbard Tree $H$ of $P$ so that it satisfies two conditions. 
First we assume that for every 
Fatou vertex $v_i$ in the Hubbard Tree, 
$H$ includes all internal rays zero ${\cal R}(i)$ incident at $v_i$. 
Also we assume that if $v$ is any vertex adjacent to $v_i$, then necessarily $v$ belongs to the boundary of the Fatou component with center $v_i$. 
In this way, given a Fatou point $v_i$, 
the set of vertices adjacent to $v_i$ can be indexed by a rational modulo $1$. 
Furthermore, this rational corresponds to the argument of the internal ray which lands at such point. 
Denote by $\Theta_i$ the set of these arguments. 
\smallskip
The second step is to construct the Hubbard Forest corresponding to $f$. 
For every $v_i \in |S|$ and $\theta \in \Theta_i$, 
we make sure that the tree $T_i$ corresponding to $v_i$ includes the orbit of the vertex $v$ 
at which the ray of argument $\theta$ lands. 
Furthermore, by including extra Julia vertices (and/or by taking homogeneous coverings) 
we may assume that $T_i-\{v\}$ has the same number of connected components as $J(f) \cap {\bf C_i} - \{v\}$. 
\smallskip
The last step is to modify the Hubbard Tree $H$ as follows. 
For every $v_i \in |S|$, 
we remove $v_i$ and all edges adjacence to $v_i$ from $H$. 
Note that every edge adjacent to $v_i$ in $H$ has associated a vertex $v$ and an argument $\theta$. 
We glue again this vertex to the tree $T_i$ between the edges  at which the ray with argument $\theta$ lands. 
Clearly the resulting tree can be made an expanding Hubbard Tree after normalization of the angles at these points. 
\medskip
{\bf 6.4 Example.} 
Consider the degree three hyperbolic polynomial 
$P(z)=z^3-{3 \over 2}z$; 
for which the critical points are changed under iteration. 
As each cycle of Fatou components has degree 4, 
there are three possible choices of `zero internal rays'. 
\medskip
\centerline{\psfig{figure=forest.fi.01.ps,height=1.6in}}
\smallskip
\vfill\eject
Let $S$ be the mapping schema of $P$. 
We consider the Hyperbolic map $f \in {\cal P}^S$ with Hubbard forest 
\vskip 1.2in 
We extend this Hubbard forest so that it contains all periodic vertices of return period 1 and rotation number zero. 
We proceed now to tune $P$ with $f$. 
We have three different possible ways to do this, 
according to the choice of zero internal rays as $[x_0,p_{10}]$,  $[x_0,p_{20}]$, $[x_0,p_{30}]$ respectively. 
The solutions with the respective remarks are shown below. 
\medskip

{\it Figure 2. The internal zero rays are $[x_0,p_{10}]$ and $[x_1,p_{10}]$.}

\end

We will need a lemma which are not stated in [P1]. 
In what follows $M$ will denote any invariant set containing the critical points of $P$. 
\medskip
{\bf 3.5. Lemma.} 
{\it Let $T^{*}_{P,M}$ be the family whose elements are closures of the components of $T_{P,M}-V_{P,M}$. 
Let $z,z' \in J(P)$ be such that for every $n \ge 0$ the iterates $P^{\circ n}(z),P^{\circ n}(z')$ belong to the same $T_n \in T^*$; then $z=z'$.} 
\medskip
{\bf Proof.} 
Suppose $z \ne z'$. 
Let $[z,z']$ be the regulated arc joining $z$ and $z'$. 
By hypothesis for every point $z_1 \in [z,z']$ different from $z,z'$, 
the successive iterates $P(z_1),P^{\circ 2}(z_1),\dots$ are all disjoint from $V$. 
In particular it follows that $J(P)-\{z_1\}$ has exactly two components. 
Next we choose a small neighborhood $U$ of $z_1$. 
The successive iterates $P(U),P^{\circ 2}(U),\dots$ must eventually intersect $V$. 
If $U$ was chosen small enough, 
it follows (because $J(P)$ is locally connected) 
that there is a $z_2 \in [z,z']$ (different from $z,z'$) which eventually maps to $V$, which is a contradiction.
\endofproof

\medskip
In section 6 as an application of the results obtains so far we 
generalize the Douady-Hubbard 
construction of `modulation' or `tuning' 
for quadratic maps.

{\bf Proof of Lemma 4.9.} 
Fix an integer $k\ge 1$ and denote by $\Lambda$ the set  of arguments $\theta$ for 
which the argument $d(u)\theta$ 
is periodic under the standard $m$-fold multiplication.  
Because of Proposition 4.2 we may assume\wlog 
that the landing point of every ray $R_\theta$ with $\theta \in \Lambda$ is in $V_u$. 
Now let $p$ be the landing point of a ray of period $k$ (for $P_u$). 
Let $\theta_1,\dots,\theta_r$ be the arguments of the rays which land at $p$. 
When we walk both trees counterclockwise starting at the $0$ rays $R_0$ and $\alpha(R_0)$, 
and compare the relative position of the rays with argument in $\Lambda$, 
we conclude that 
the rays $\alpha(R_{\theta_1}),\dots,\alpha(R_{\theta_1})$ land at $\alpha(p)$; 
and furthermore, only these rays land there. 
\smallskip
An analogous argument shows that 
whenever $q$ is the landing point of the rational rays $R_{\theta_1},\dots,R_{\theta_r}$, 
the corresponding rays  
$\alpha(R_{\theta_1}),\dots,\alpha(R_{\theta_r})$ land at the corresponding point $\alpha(q)$. 
Furthermore, only these rays land at $\alpha(q)$. 
\endofproof
$S_f=(|S_f|,F_f,w_f)$ (compare [M] and Appendix A); 
where $|S_f|$ is a set of vertices in precise correspondence with the connected components of $M$ (i.e, with each copy of the complex plane). 
The functions $F_f$ and $w_f$ which map $|S_f|$ respectively to itself and the non negative integers, 
are defined as follows. 
Take $u \in |S_f|$, then $f$ restricted to the copy ${\bf C}_u$ is a polynomial of degree $w_f(u)+1$ taking values in the copy ${\bf C}_{F_f(u)}$.